\documentclass{article}

\usepackage[utf8]{inputenc}
\usepackage{subfigure}
\usepackage{amsmath}
\usepackage{amssymb}
\usepackage{amsthm}
\usepackage[pdftex]{hyperref}
\usepackage{tikz}
\usepackage{caption}
\usepackage[round]{natbib}
\usepackage{fancyhdr}
\usepackage{amsfonts}
\usepackage{times}
\usepackage{ifpdf}
\usepackage{latexsym}
\usepackage{graphicx}
\usepackage{enumerate}

\newtheorem{definition}{Definition}[section]

\newtheorem{theorem}[definition]{Theorem}
\newtheorem{lemma}[definition]{Lemma}
\newtheorem{corollary}[definition]{Corollary}

\newlength{\taille} \makeatletter
\def\qed{%
  \ifmmode\vrule width .5\baselineskip height 0pt depth .5\baselineskip%
  \else{%
    \unskip\nobreak\hfil%
    \setlength{\taille}{\f@size\p@}%
    \penalty50\hskip1em\null\nobreak\hfil\vrule width .5\taille height
    0pt depth .5\taille
    \parfillskip=0pt\finalhyphendemerits=0\endgraf}%
  \fi} \makeatother

\newlength{\taillepreuve}
\newenvironment{demo}{%
  \setbox123=\hbox{Proof:}%
  \taillepreuve=\wd123%

  \vskip-\lastskip\nobreak\medskip\par\noindent\box123\list{}{\leftmargin
    .5\taillepreuve}\parindent=1em\item} {\qed\endlist\bigskip}

\DeclareMathOperator{\Domain}{Domain}
\DeclareMathOperator{\OrderFunction}{OrderFunction}

\DeclareMathOperator{\Height}{Height}
\DeclareMathOperator{\Level}{Level}

\DeclareMathOperator{\Inv}{Inv}

\DeclareMathOperator{\Sup}{Supremum}
\DeclareMathOperator{\Infimum}{Infimum}

\DeclareMathOperator{\DisjointCopy}{DisjointCopy}
\DeclareMathOperator{\Gap}{Gap}
\DeclareMathOperator{\InitialSection}{InitialSection}

\begin{document}



\author{Laurent Lyaudet\footnote{\url{https://lyaudet.eu/laurent/}, laurent.lyaudet@gmail.com}}
\title{On level-induced suborders}

\maketitle
\begin{abstract}
In this article, we characterize orders that are level-induced suborders anytime they are induced suborders of a superorder.
We also characterize orders that are consecutive level-induced suborders anytime they are level-induced suborders of a superorder.
Thus characterizing orders that are consecutive level-induced suborders anytime they are induced suborders of a superorder.
\end{abstract}

Current version : 2020/03/29


\section{Introduction}
\label{section:introduction}

Apologies: We do science as a hobby, it is not our daily job and there is an impact on the quality of the bibliography.
For an unpublished work we did in 2015, we started doing bibliographic search during 9 months,
but all the gathered references were lost when a hacker erased all our files on our laptop.
Since then, we chose to publish our ideas on arXiv and correct the bibliography afterwards.
For this article, we found no prior work defining kinds of induced suborders with constraints
on their levels relatively to those of the superorder.
Our search was in English and French scientific literature,
and since the topic of order theory is ancient and vast,
we may have missed early references in other languages.
If you do know an early reference, please be kind enough to email/correct us.
This is now the fifth version on arXiv, we are sorry for the errors that we published in the preceding four versions;
we hope there is none in this version; at least, there is some (slow) progress.

This article study Open problem 5.16 in \cite{DBLP:journals/corr/abs-1903-02028}.
``Characterize finite orders that are induced suborders of any
well-founded order if and only if they are (consecutive) level-induced suborders of this
well-founded order. Examples: chains, antichains of size 1 and 2. Counter-examples:
antichains of size at least 3.''

Section \ref{section:definitions_and_notations} contains most of the definitions and notations used in this article.
In section \ref{section:Orders_that_are_always_level-induced_suborders},
we characterize orders that cannot be induced suborders without being level-induced suborders.
Section \ref{section:Orders_that_are_naturally_consecutive_level-induced_suborders} characterize orders
that cannot be level-induced suborders without being consecutive level-induced suborders.
In section \ref{section:Algorithms_to_find_ali_induced_suborders_and_nacli_level-induced_suborders},
we give algorithms to find (level-)induced suborders of the previously defined classes.

\section{Definitions and notations}
\label{section:definitions_and_notations}

Throughout this article, we use the following definitions and notations.
\begin{math}O\end{math} will be reserved for asymptotic growth of functions.
Thus \begin{math}P\end{math} denotes an order (it may be either a partial, or a total/linear order),
in particular \begin{math}P^{0,1}\end{math} denotes the binary total order where \begin{math}0 < 1\end{math}.
We denote \begin{math}\Domain(P)\end{math}, the domain of the order \begin{math}P\end{math}
(for example, \begin{math}\Domain(P^{0,1}) = \{0,1\}\end{math}).
We write \begin{math}x < y\end{math}, and \begin{math}x > y\end{math} as usual to express the order between two elements;
we also write \begin{math}x \sim y\end{math} when two elements are incomparable in the partial order considered.
We denote \begin{math}\OrderFunction(P)\end{math}, the order function of the order \begin{math}P\end{math}
defined from \begin{math}\Domain(P)^2\end{math} to \begin{math}\{=,\sim,<,>\}\end{math}
(for example, \begin{math}\OrderFunction(P^{0,1}) = \{((0,0),=),~((0,1),<),~((1,0),>),~((1,1),=)\}\end{math}).

We denote \begin{math}\Inv(P)\end{math}, the inverse/reverse order of \begin{math}P\end{math};
for example, \begin{math}\Inv(P^{0,1}) = P^{1,0}\end{math} is the order on 0 and 1
where \begin{math}1 < 0\end{math}.

\begin{definition}[Maximum chain, height]
Let \begin{math}P\end{math} be an order,
a chain of \begin{math}P\end{math} is \emph{maximum}
if it is maximal and no other chain of \begin{math}P\end{math} has greater cardinality.
The cardinal of a maximum chain is the \emph{height} of \begin{math}P\end{math},
denoted \begin{math}\Height(P)\end{math}.
When \begin{math}P\end{math} is well-founded\footnote{
Some authors also say Noetherian. In both cases, it means that there is no strictly decreasing infinite sequence.
},
we redefine a maximum chain to be one such that the corresponding ordinal is maximum;
and we redefine its height to be the ordinal corresponding to its maximum chains.
Thus in this case \begin{math}\Height(P)\end{math} denotes an ordinal.
\end{definition}
Note that an infinite order may have no maximum chain, but it always have at least one maximal chain.
When there is no maximum chain, \begin{math}\Height(P)\end{math} is defined as
the supremum cardinal/ordinal of the cardinals/ordinals corresponding to maximal chains.

In a well-founded order \begin{math}P\end{math},
the level decomposition of \begin{math}P\end{math}
is the function \begin{math}\Level_{P}: \Domain(P) \rightarrow \Height(P)\end{math}
(\begin{math}\Height(P)\end{math} is an arbitrary ordinal.)
such that \begin{math}\forall x \in \Domain(P), \Level_{P}(x) = \Sup(\Level_{P}(y) + 1 \text{ such that } y < x, y \in \Domain(P))\end{math}.
(Of course, this supremum is 0 if no element is below \begin{math}x\end{math}.)
We define the \emph{level-width} of P as the supremum of the cardinals of the levels of P.
Given two elements \begin{math}x, y \in \Domain(P), \Gap_{P}(x,y) = \Gap_{P}(y,x) = \Sup(\Level_{P}(x), \Level_{P}(y)) - \Infimum(\Level_{P}(x), \Level_{P}(y))\end{math}.
The gap between two elements is clearly 0 if and only if these two elements belong to the same level.
Note that this ordinal gap may not correspond to an actual well-founded chain in \begin{math}P\end{math} between the two elements,
in particular they may have an arbitrary large gap and be incomparable.
When there is more than two elements, the gap of a set of elements is the supremum of the gaps of the pairs.

We consider the following kinds of suborders:
\begin{itemize}
\item An induced suborder \begin{math}P'\end{math} of an order \begin{math}P\end{math}
      is such that \begin{math}\Domain(P') \subseteq \Domain(P)\end{math},
      and \begin{math}\forall x,y \in \Domain(P'), \OrderFunction(P')(x,y) = \OrderFunction(P)(x,y)\end{math}.
      Let \begin{math}X \subseteq \Domain(P)\end{math},
      \begin{math}P[X]\end{math} denotes the suborder of \begin{math}P\end{math} induced by \begin{math}X\end{math}.
\item A \emph{level-induced suborder} \begin{math}P'\end{math} of a well-founded order \begin{math}P\end{math}
      is such that \begin{math}\Domain(P') \subseteq \Domain(P)\end{math},
      \begin{math}\forall x,y \in \Domain(P'), \OrderFunction(P')(x,y) = \OrderFunction(P)(x,y)\end{math},
      and \begin{math}\forall x,y \in \Domain(P'), \Level_{P'}(x) = \Level_{P'}(y) \Leftrightarrow \Level_{P}(x) = \Level_{P}(y)\end{math}.
      (Note that we could also define two other kinds of level-induced suborders with
      \begin{math}\forall x,y \in \Domain(P'), \Level_{P'}(x) = \Level_{P'}(y) \Rightarrow \Level_{P}(x) = \Level_{P}(y)\end{math}
      \begin{math}\forall x,y \in \Domain(P'), \Level_{P'}(x) = \Level_{P'}(y) \Leftarrow \Level_{P}(x) = \Level_{P}(y)\end{math}.
      There is a simple proof by transfinite induction on the levels of \begin{math}P'\end{math} showing that
      \begin{math}\forall x,y \in \Domain(P'), \Level_{P'}(x) = \Level_{P'}(y) \Rightarrow \Level_{P}(x) = \Level_{P}(y)\end{math}
      implies \begin{math}\forall x,y \in \Domain(P'), \Level_{P'}(x) = \Level_{P'}(y) \Leftarrow \Level_{P}(x) = \Level_{P}(y)\end{math}.
      Moreover, the same proof shows that
      \begin{math}\forall x,y \in \Domain(P'), \Level_{P'}(x) < \Level_{P'}(y) \Leftrightarrow \Level_{P}(x) < \Level_{P}(y)\end{math}.
      Thus only two kinds of level-induced suborders exists
      \begin{math}\forall x,y \in \Domain(P'), \Level_{P'}(x) = \Level_{P'}(y) (\Leftrightarrow \text{ or } \Rightarrow) \Level_{P}(x) = \Level_{P}(y)\end{math},
      and \begin{math}\forall x,y \in \Domain(P'), \Level_{P'}(x) = \Level_{P'}(y) \Leftarrow \Level_{P}(x) = \Level_{P}(y)\end{math}.  )
\item A \emph{consecutive level-induced suborder} \begin{math}P'\end{math} of a well-founded order \begin{math}P\end{math}
      is such that \begin{math}\Domain(P') \subseteq \Domain(P)\end{math},
      \begin{math}\forall x,y \in \Domain(P'), \OrderFunction(P')(x,y) = \OrderFunction(P)(x,y)\end{math},
      \begin{math}\forall x,y \in \Domain(P'), \Gap_{P'}(x,y) = \Gap_{P}(x,y)\end{math}.
      (In the finite case, the equality of the gaps may be replaced by the following conditions:
      \begin{math}\forall x,y \in \Domain(P'), \Level_{P'}(x) = \Level_{P'}(y) \Leftrightarrow \Level_{P}(x) = \Level_{P}(y)\end{math},
      and \begin{math}\forall x,y \in \Domain(P'), \Level_{P'}(x) + 1 = \Level_{P'}(y) \Leftrightarrow \Level_{P}(x) + 1 = \Level_{P}(y)\end{math}.)
\end{itemize}

\section{Orders that are always level-induced suborders}
\label{section:Orders_that_are_always_level-induced_suborders}

In this section, we assume that a given well-founded order \begin{math}P'\end{math}
is an induced suborder of a well-founded order \begin{math}P\end{math}.
We study necessary and sufficient conditions on \begin{math}P'\end{math}
to have that \begin{math}P'\end{math}
is a level-induced suborder of \begin{math}P\end{math}.
\begin{definition}[ali orders]
An \emph{ali order} \begin{math}P'\end{math} is a well-founded order
such that whenever \begin{math}P'\end{math} is isomorphic to an induced suborder of a well-founded order \begin{math}P\end{math},
then \begin{math}P'\end{math} is also isomorphic to a level-induced suborder of \begin{math}P\end{math}.
\end{definition}

Recall that an initial section \begin{math}I\end{math} of an order \begin{math}P\end{math}
is a subset of his domain closed by taking smaller elements: \begin{math}\forall x \in \Domain(P), \forall y \in I, x < y \Rightarrow x \in I\end{math}.
Given a subdomain \begin{math}Y\end{math} of \begin{math}P\end{math},
\begin{math}\InitialSection(P,Y) = \{ x \in \Domain(P), \exists y \in Y, x \leq y\}\end{math}
is the initial section generated by \begin{math}Y\end{math}.
(\begin{math}I = \InitialSection(P,I)\end{math}
and \begin{math}P[I]\end{math} denotes the suborder of \begin{math}P\end{math} induced by \begin{math}I\end{math}.)

\begin{lemma}
Let \begin{math}P'\end{math} be a well-founded order.
\begin{itemize}
\item If we have a level-induced suborder isomorphic to \begin{math}P'\end{math}
      in the suborder induced by an initial section of some well-founded order \begin{math}P\end{math},
      then this level-induced suborder is also level-induced in \begin{math}P\end{math}.
\item If \begin{math}P'\end{math} is an ali order and it is isomorphic to an induced suborder \begin{math}P''\end{math}
      of the well-founded order \begin{math}P\end{math},
      it is isomorphic to a level-induced suborder in the restriction of \begin{math}P\end{math}
      to the initial section generated by \begin{math}P''\end{math}.
\end{itemize}
\end{lemma}
\begin{demo}
The first assertion is trivially true because
\begin{math}\forall x \in \InitialSection(P,X)\end{math}, \newline\begin{math}\Level_{P[\InitialSection(P,X)]}(x) = \Level_{P}(x)\end{math}.

The second assertion follows from it and the fact that \begin{math}P'\end{math} is an ali order
(``relatively to \begin{math}P[\InitialSection(P,\Domain(P''))]\end{math}'').
\end{demo}

\begin{corollary}[ali orders revisited]
An \emph{ali order} \begin{math}P'\end{math} is a well-founded order
such that whenever \begin{math}P'\end{math} is isomorphic to an induced suborder \begin{math}P''\end{math} of a well-founded order \begin{math}P\end{math},
then \begin{math}P'\end{math} is also isomorphic to a level-induced suborder of \begin{math}P[\InitialSection(P,\Domain(P''))]\end{math}.
\end{corollary}

We now observe that :
\begin{lemma}
Any well-founded order \begin{math}P'\end{math} is an induced suborder
of a well-founded order \begin{math}P\end{math} of level-width 2.
Moreover, \begin{math}P\end{math} has no level-induced suborder isomorphic
to \begin{math}O_{obs1} = (\{a,b,c,d\}, \{a < b, c < d\}) \equiv \Inv(O_{obs1})\end{math}.
\end{lemma}
\begin{demo}
We use a well-founded chain to lift each element of \begin{math}\Domain(P')\end{math} to a separate level.
By Zermelo's axiom, there is a bijection \begin{math}f\end{math}
between some ordinal \begin{math}\alpha\end{math} and \begin{math}\Domain(P')\end{math},
such that \begin{math}\forall x,y \in \Domain(P'), \Level_{P'}(x) < \Level_{P'}(y) \Rightarrow f(x) < f(y)\end{math}.
Assume, without loss of generality, that \begin{math}\Domain(P') \cap \alpha = \emptyset\end{math}.

Let \begin{math}\Domain(P) = \Domain(P') \sqcup \alpha\end{math}:
\begin{itemize}
\item \begin{math}\forall x,y \in \Domain(P'), \OrderFunction(P)(x,y) = \OrderFunction(P')(x,y)\end{math},
\item \begin{math}\forall x,y \in \alpha, \OrderFunction(P)(x,y) = \OrderFunction(\alpha)(x,y)\end{math},
\item \begin{math}\forall x \in \alpha, \forall y \in \Domain(P'), \OrderFunction(P)(x,y) = \text{ `}<\text{' if } x < f^{-1}(y), \text{ `}\sim\text{' otherwise }\end{math}.
\end{itemize}
Clearly, \begin{math}P\end{math} has all the claimed properties.
\end{demo}

\begin{corollary}
An ali order has level-width at most 2 and no level-induced suborder isomorphic
to \begin{math}O_{obs1} = (\{a,b,c,d\}, \{a < b, c < d\}) \equiv \Inv(O_{obs1})\end{math}.
\end{corollary}

\begin{lemma}
No ali order has a level of size 2 except maybe the first.
\end{lemma}
\begin{demo}
Assume for a contradiction that \begin{math}P'\end{math} is an ali order with two elements
\begin{math}x,y\end{math} such that \begin{math}\Level_{P'}(x) = \Level_{P'}(y) > 0\end{math}.
Take \begin{math}x,y\end{math} such that their level is minimum.
By the previous corollary, we must have an element \begin{math}z, z < x, z < y\end{math}.
(By transitivity, it is trivial to see that such a \begin{math}z\end{math} exists in all previous levels,
since only the first level may have two elements.)
Thus we have a level-induced suborder isomorphic to \begin{math}O_{obs2} = (\{a,b,c\}, \{a < b, a < c\})\end{math}.

We now show how to remove all such level-induced suborders for any well-founded order.
Again by Zermelo's axiom, there is a bijection \begin{math}f\end{math}
between some ordinal \begin{math}\alpha\end{math} and \begin{math}\Domain(P')\end{math},
such that \begin{math}\forall x,y \in \Domain(P'), \Level_{P'}(x) < \Level_{P'}(y) \Rightarrow f(x) < f(y)\end{math}.
This time, we add a distinct chain for each element of \begin{math}\Domain(P')\end{math}.
Let \begin{math}\DisjointCopy(i)\end{math} be a chain isomorphic to the ordinal \begin{math}i\end{math}
such that its elements are assumed to be distinct from all other elements considered in the following formula:
\begin{math}\Domain(P) = \Domain(P') \sqcup (\bigsqcup_{i \in \alpha} \DisjointCopy(i))\end{math}.
\begin{itemize}
\item \begin{math}\forall x,y \in \Domain(P'), \OrderFunction(P)(x,y) = \OrderFunction(P')(x,y)\end{math},
\item \begin{math}\forall x,y \in \DisjointCopy(i), \OrderFunction(P)(x,y) = \OrderFunction(\DisjointCopy(i))(x,y)\end{math},
\item \begin{math}\forall x \in \DisjointCopy(i), \forall y \in \DisjointCopy(j), \OrderFunction(P)(x,y) = \text{ `}\sim\text{'}\end{math},
\item \begin{math}\forall x \in \DisjointCopy(i), \forall y \in \Domain(P'),
        \OrderFunction(P)(x,y) = \text{ `}<\text{' if } i = f^{-1}(y) \text{ or } f(i) < y (\OrderFunction(P')(f(i),y) \in \{=, <\}), \text{ `}\sim\text{' otherwise }\end{math}.
\end{itemize}
Clearly, each element \begin{math}f(i)\end{math} of \begin{math}\Domain(P')\end{math} is now on a distinct level,
since \begin{math}\DisjointCopy(i)\end{math} is a longest chain below it.
Moreover, if some element in \begin{math}\DisjointCopy(i)\end{math} is less than two elements on the same level,
then clearly, \begin{math}f(i)\end{math} must be less than these two elements,
and this is impossible since \begin{math}f(i)\end{math} may only be less than elements in \begin{math}\Domain(P')\end{math},
that are now scattered.
\end{demo}

\begin{theorem}
An ali order is either
\begin{itemize}
\item a well-founded total order,
\item an antichain of size 2,
\item the disjoint union of a well-founded chain and an incomparable element,
      where the well-founded chain has height 2 or is isomorphic to a regular cardinal/ordinal or its successor,
\item the order composition of two incomparable elements and a well-founded chain
      (we call this case a ``(1,1)-based chain''),
\item or the order composition of
  \begin{itemize}
  \item the disjoint union of a well-founded chain and an incomparable element,
        where the well-founded chain has height 2 (we call this case a ``(2,1)-based chain'') or is isomorphic to a regular cardinal/ordinal or its successor,
  \item and a well-founded chain.
  \end{itemize}
\end{itemize}

\end{theorem}
\begin{demo}
It is trivial to see that a well-founded total order is an ali order.

A well-founded chain and an incomparable element may or may not form an ali order.
\begin{itemize}
\item Clearly an antichain of size 2 is an ali order.
\item A chain of height 2 and an incomparable element is an ali order:
  Indeed, consider \begin{math}x,y,z\end{math} with \begin{math}x < y, x \sim z, y \sim z\end{math}.
  \begin{itemize}
  \item If \begin{math}\Level_{P}(x) = \Level_{P}(z)\end{math}, there is nothing to do.
  \item If \begin{math}\Level_{P}(x) < \Level_{P}(z)\end{math},
        then there is another element \begin{math}t\end{math} such that \begin{math}\Level_{P}(x) = \Level_{P}(t)\end{math},
        and \begin{math}t < z\end{math}
        (since there is an element in each level below \begin{math}z\end{math},
        such that \begin{math}z\end{math} is more than it, and \begin{math}x \sim z\end{math}).
        Clearly, \begin{math}t,z,x\end{math} give a level-induced suborder isomorphic to a chain of height 2 and an incomparable element.
  \item If \begin{math}\Level_{P}(x) > \Level_{P}(z)\end{math},
        then there is another element \begin{math}t\end{math} such that \begin{math}\Level_{P}(z) = \Level_{P}(t)\end{math},
        and \begin{math}t < x\end{math}
        (since there is an element in each level below \begin{math}x\end{math},
        such that \begin{math}x\end{math} is more than it, and \begin{math}x \sim z\end{math}).
        Clearly, \begin{math}t,y,z\end{math} (or \begin{math}t,x,z\end{math}) give a level-induced suborder isomorphic to a chain of height 2 and an incomparable element.
        \emph{Note that the same argument applies to well-founded chains of any height
        with a lowest element \begin{math}x\end{math} and an incomparable element \begin{math}z\end{math}.}
  \end{itemize}
\item A well-founded chain corresponding to a successor of a successor ordinal \begin{math}\alpha + 2\end{math}
      more than 2 (\begin{math}\alpha > 0\end{math}) and an incomparable element is not an ali order.
      (Such a chain is ended by a chain of height 2 on two consecutive levels \begin{math}\alpha\end{math} and \begin{math}\alpha + 1\end{math}.)
      Indeed, consider the order \begin{math}P\end{math}
      with \begin{math}\Domain(P) = \DisjointCopy(\alpha + 2) \sqcup \DisjointCopy(\alpha + 1)\end{math},
      such that:
  \begin{itemize}
  \item \begin{math}\forall x,y \in \DisjointCopy(\alpha + 2), \OrderFunction(P)(x,y) = \OrderFunction(\DisjointCopy(\alpha + 2))(x,y)\end{math},
  \item \begin{math}\forall x,y \in \DisjointCopy(\alpha + 1), \OrderFunction(P)(x,y) = \OrderFunction(\DisjointCopy(\alpha + 1))(x,y)\end{math},
  \item \begin{math}\forall x \in \DisjointCopy(\alpha + 1), \forall y \in \DisjointCopy(\alpha + 2),
        \OrderFunction(P)(x,y) = \text{ `}<\text{' if }x\text{ corresponds to an ordinal less than }y
        \text{ and x does not correspond to }\alpha, \text{ `}\sim\text{' otherwise }\end{math}.
  \end{itemize}
  Any chain in \begin{math}P\end{math} isomorphic to \begin{math}\alpha + 2\end{math} must have its greatest element
  to be the greatest element corresponding to \begin{math}\alpha + 1\end{math} in \begin{math}\DisjointCopy(\alpha + 2)\end{math}.
  Then clearly the only element incomparable with it is the greatest element corresponding to \begin{math}\alpha\end{math} in \begin{math}\DisjointCopy(\alpha + 1)\end{math}.
  Thus, since \begin{math}\alpha > 0\end{math} it cannot yield a level-induced suborder.

\item A well-founded chain corresponding to (a successor of) a limit ordinal \begin{math}\alpha\end{math}
      that is singular (not a regular cardinal/ordinal)
      and an incomparable element is not an ali order.
      (Such a chain is isomorphic to \begin{math}\gamma = \alpha\end{math} (resp. \begin{math}\gamma = \alpha + 1\end{math}).)
      Since \begin{math}\alpha\end{math} is not a regular cardinal,
      let \begin{math}\beta + 1 < \gamma\end{math} be a successor ordinal such that
      \begin{math}\alpha \setminus \beta\end{math} is not isomorphic to \begin{math}\alpha\end{math}.
      Now, consider the order \begin{math}P\end{math}
      with \begin{math}\Domain(P) = \DisjointCopy(\gamma) \sqcup \DisjointCopy(\beta + 1)\end{math},
      such that:
  \begin{itemize}
  \item \begin{math}\forall x,y \in \DisjointCopy(\gamma), \OrderFunction(P)(x,y) = \OrderFunction(\DisjointCopy(\gamma))(x,y)\end{math},
  \item \begin{math}\forall x,y \in \DisjointCopy(\beta + 1), \OrderFunction(P)(x,y) = \OrderFunction(\DisjointCopy(\beta + 1))(x,y)\end{math},
  \item \begin{math}\forall x \in \DisjointCopy(\beta + 1), \forall y \in \DisjointCopy(\gamma),
        \OrderFunction(P)(x,y) = \text{ `}<\text{' if }x\text{ corresponds to an ordinal less than }y
        \text{ and x does not correspond to }\beta, \text{ `}\sim\text{' otherwise }\end{math}.
  \end{itemize}
  Since \begin{math}\beta + 1 < \gamma\end{math}, any chain in \begin{math}P\end{math} isomorphic to
  \begin{math}\gamma\end{math} must have its final segment in \begin{math}\DisjointCopy(\gamma)\end{math}.
  Then clearly the only element incomparable with it is the greatest element corresponding to \begin{math}\beta\end{math} in \begin{math}\DisjointCopy(\beta + 1)\end{math}.
  Thus, since \begin{math}\beta > 0\end{math},
  and \begin{math}\alpha \setminus \beta\end{math} is not isomorphic to \begin{math}\alpha\end{math},
  it cannot yield a level-induced suborder.

\item A well-founded chain corresponding to (a successor of) a regular limit ordinal/cardinal \begin{math}\alpha\end{math}
      and an incomparable element is an ali order.
      (Such a chain is isomorphic to \begin{math}\gamma = \alpha\end{math} (resp. \begin{math}\gamma = \alpha + 1\end{math}).)
      Let \begin{math}(x_i)_{i \in \gamma}\end{math} and an incomparable element \begin{math}z\end{math}
      form such an induced suborder in some partial order \begin{math}P\end{math}.
      From what we noted for the case of a chain of height 2 and an incomparable element,
      we just need to consider the case where \begin{math}\Level_{P}(x_0) < \Level_{P}(z)\end{math}.
      Clearly, there is a chain \begin{math}(y_j)_{j \in \Level_{P}(z)}\end{math} below \begin{math}z\end{math}
      intersecting all levels below \begin{math}\Level_{P}(z)\end{math},
      such that, by transitivity, \begin{math}x_i \not< y_j, i \in \gamma, j \in \Level_{P}(z)\end{math}.
      \begin{itemize}
      \item If \begin{math}\Level_{P}(z) \geq \Sup(\Level_{P}(x_i), i \in \gamma)\end{math},
            then \begin{math}(y_j)_{j \in \Level_{P}(z), j \geq \Level_{P}(x_0)}\end{math}
            contains a subchain isomorphic to \begin{math}\gamma\end{math},
            and \begin{math}x_0 \sim y_j, j \in \Level_{P}(z), j \geq \Level_{P}(x_0)\end{math}.
            Taking the subchain starting on the same level than \begin{math}x_0\end{math}
            (made of the elements \begin{math}(y_j)_{j \in \Level_{P}(z), j \geq \Level_{P}(x_0), j \in \{\Level(x_i), i \in \gamma\}}\end{math}),
            together with \begin{math}z\end{math} on top of this subchain (if needed) and \begin{math}x_0\end{math},
            we obtain the sought level-induced suborder.
      \item If \begin{math}\Level_{P}(x_0) < \Level_{P}(z) < \Sup(\Level_{P}(x_i), i \in \gamma)\end{math}, it is more complicated.
            If \begin{math}(x_i)_{i \in \gamma, \Level_{P}(x_i) \geq \Level_{P}(z)}\end{math}
            is reduced to the singleton \begin{math}x_{\alpha}\end{math} (hence \begin{math}\gamma = \alpha + 1\end{math}),
            then we are in a situation equivalent to the previous case \begin{math}\Level_{P}(z) \geq \Sup(\Level_{P}(x_i), i \in \gamma)\end{math},
            but we should replace it by \begin{math}\Level_{P}(z) > \Level_{P}(x_i), \forall i \in \alpha\end{math};
            it is clear in that case that \begin{math}(y_j)_{j \in \Level_{P}(z), j \geq \Level_{P}(x_0)}\end{math}
            contains a subchain isomorphic to \begin{math}\alpha\end{math} starting on the same level than \begin{math}x_0\end{math};
            together with \begin{math}z\end{math} on top of this subchain and \begin{math}x_0\end{math},
            it yields the sought level-induced suborder.

            Assume that no element \begin{math}x_i, i \in \gamma\end{math} belongs to \begin{math}\Level_{P}(z)\end{math}.
            There may be no element in \begin{math}\Level_{P}(z)\end{math} ordered with all elements \begin{math}x_i, i \in \gamma\end{math};
            nevertheless, there is a lowest element \begin{math}x_k\end{math} such that \begin{math}\Level_{P}(x_k) > \Level_{P}(z)\end{math};
            and there is an element \begin{math}x'_k < x_k\end{math} such that \begin{math}\Level_{P}(x'_k) = \Level_{P}(z)\end{math}.
            If there is an element \begin{math}x_i, i \in \gamma\end{math} that belongs to \begin{math}\Level_{P}(z)\end{math},
            we also name this element \begin{math}x'_k\end{math}.

            Clearly, \begin{math}\{x'_k\} \cup (x_i)_{i \in \gamma, \Level_{P}(x_i) \geq \Level_{P}(z)}\end{math}
            is cofinal in a chain isomorphic to \begin{math}\gamma\end{math},
            and all elements of the chain are incomparable with \begin{math}z\end{math}.
            Moreover, even in the case \begin{math}\gamma = \alpha + 1\end{math},
            we can remove the element \begin{math}x_{\alpha}\end{math} of this chain,
            and still obtain a chain cofinal in a chain isomorphic to \begin{math}\alpha\end{math},
            since we already studied the case where \begin{math}(x_i)_{i \in \gamma, \Level_{P}(x_i) \geq \Level_{P}(z)}\end{math}
            is reduced to the singleton \begin{math}x_{\alpha}\end{math}.
            In both cases, by regularity of \begin{math}\alpha\end{math},
            the subchain \begin{math}\{x'_k\} \cup (x_i)_{i \in \gamma, \Level_{P}(x_i) \geq \Level_{P}(z)}\end{math}
            is isomorphic to \begin{math}\gamma\end{math}.
            Thus, together with the incomparable element \begin{math}z\end{math},
            it yields the sought level-induced suborder.
      \end{itemize}
\end{itemize}

The only case left to study is then when the second element on level 0, \begin{math}y\end{math},
is less than some element of the chain \begin{math}(x_i)_{i \in \gamma}\end{math}.
Let \begin{math}(x_i)_{i \in \alpha}\end{math} be the initial segment of elements
that are incomparable with \begin{math}y\end{math}.

Let \begin{math}P\end{math} be an order containing the chain \begin{math}(x_i)_{i \in \gamma}\end{math}
and the element \begin{math}y\end{math} with appropriate order relationship.
If \begin{math}(x_i)_{i \in \alpha}\end{math} and \begin{math}y\end{math} form an ali order,
then clearly by transitivity all elements in \begin{math}P[\InitialSection(P,\{x_i, i \in \alpha\} \cup \{y\})]\end{math}
will be less than elements \begin{math}(x_i)_{i \in \gamma \setminus \alpha}\end{math}.
Thus, a level-induced suborder isomorphic to \begin{math}(x_i)_{i \in \alpha}\end{math} and \begin{math}y\end{math}
in \begin{math}P[\InitialSection(P,\{x_i, i \in \alpha\} \cup \{y\})]\end{math},
together with the chain \begin{math}(x_i)_{i \in \gamma \setminus \alpha}\end{math},
will immediately yield a level-induced suborder isomorphic to \begin{math}(x_i)_{i \in \gamma}\end{math} and \begin{math}y\end{math}.

If \begin{math}(x_i)_{i \in \alpha}\end{math} and \begin{math}y\end{math} does not form an ali order,
consider an order \begin{math}P_{counter}\end{math} containing \begin{math}(x_i)_{i \in \alpha}\end{math} and \begin{math}y\end{math}
with appropriate order relationship such that no level-induced suborder isomorphic to
\begin{math}(x_i)_{i \in \alpha}\end{math} and \begin{math}y\end{math} exists.
Let \begin{math}P\end{math} be the order composition of \begin{math}P_{counter}\end{math}
and the chain \begin{math}(x_i)_{i \in \gamma \setminus \alpha}\end{math}.
Clearly, any element of \begin{math}(x_i)_{i \in \gamma \setminus \alpha}\end{math} is comparable with all other elements
of \begin{math}P\end{math}.
Thus, no (level-)induced suborder of \begin{math}P\end{math} isomorphic to
\begin{math}(x_i)_{i \in \alpha}\end{math} and \begin{math}y\end{math}
may contain an element of \begin{math}(x_i)_{i \in \gamma \setminus \alpha}\end{math}.
Hence, no level-induced suborder isomorphic to
\begin{math}(x_i)_{i \in \alpha}\end{math} and \begin{math}y\end{math}
exists in \begin{math}P\end{math},
and no level-induced suborder isomorphic to
\begin{math}(x_i)_{i \in \gamma}\end{math} and \begin{math}y\end{math}
exists in \begin{math}P\end{math}.

This ends this technical proof.
\end{demo}

\begin{corollary}
An ali order is a \emph{well-founded} order without induced suborder isomorphic
to \begin{math}O_{obs1} = (\{a,b,c,d\}, \{a < b, c < d\}) \equiv \Inv(O_{obs1})\end{math},
\begin{math}O_{obs2} = (\{a,b,c\}, \{a < b, a < c\})\end{math}, or an antichain of size 3,
but infinite ali orders cannot be characterized by a set of forbidden induced suborders
(it may be a class of forbidden induced suborders but not a set, let alone a finite set,
because any chain with an incomparable element is a suborder of a regular chain with an incomparable element).
Thus, ali orders are a subclass of series parallel interval orders.

Nevertheless, finite ali orders are the finite orders without induced suborder isomorphic
to \begin{math}O_{obs1} = (\{a,b,c,d\}, \{a < b, c < d\}) \equiv \Inv(O_{obs1})\end{math},
\begin{math}O_{obs2} = (\{a,b,c\}, \{a < b, a < c\})\end{math}, an antichain of size 3,
or \begin{math}(\{a,b,c,d\}, \{a < b, a < c, b < c\})\end{math} (a chain of height 3 and an incomparable element).
Thus, finite ali orders are a subclass of the following classes:
series parallel unit interval orders,
semi-orders = unit interval orders,
1-weak orders (see \cite{DBLP:journals/dm/Trenk98}).
\end{corollary}

Thus finite ali orders can be recognized in time \begin{math}O(n+m)\end{math},
where \begin{math}n\end{math} is the number of elements
and \begin{math}m\end{math} is the number of comparability relationships,
see the articles by \cite{DBLP:conf/stoc/ValdesTL79} and \cite{DBLP:journals/dam/CrespelleP06}, for example.

(Given a modular decomposition using disjoint sum and order composition,
a simple tree-automaton can determine if it corresponds to an ali order and compute the length of the longest chain.
There are two variants of modular decomposition:
\begin{itemize}
\item the binary one where binary disjoint sum and binary order composition have exactly two subtrees/subterms,
\item the grouped one where grouped disjoint sum and grouped order composition may have more than two subtrees below,
      no two grouped disjoint sum nodes are adjacent in the decomposition,
      and no two grouped order composition nodes are adjacent in the decomposition.
\end{itemize}
The grouped variant can simplify some computations.
Below, we precise when the computation applies only to one variant.
Without loss of generality, we assume that on each node, we have a boolean bLeaf: true if the node is a leaf,
false if it is a disjoint sum node or an order composition node.
The set of states of the tree-automaton has size 12, it is the cartesian product of 3 sets of substates:
\begin{itemize}
\item a boolean value bChain which is true if and only if the order defined by the modular decomposition up to this node is a chain/total order.
      bChain of a leaf/single element is true, bChain of a disjoint sum is false,
      bChain of an order composition is a logical conjunction (an AND) of the values of bChain for the subtrees of the order composition.
\item a ternary value iHeight1\_2\_More which is 1, 2, or ``3'' if and only if the order defined by the modular decomposition up to this node has height 1, 2, or more than 2.
      iHeight1\_2\_More of a leaf/single element is 1, iHeight1\_2\_More of a disjoint sum is the maximum of the values of iHeight1\_2\_More for the subtrees of the disjoint sum,
      iHeight1\_2\_More of an order composition is the minimum of 3 and the sum of the values of iHeight1\_2\_More for the subtrees of the order composition,
      this sum with upper-bound 3 is either 2 if the order composition has only two subterms and both have iHeight1\_2\_More = 1, or 3 otherwise.
\item a boolean value bAli which is true if and only if the order defined by the modular decomposition up to this node is an ali order.
      bAli of a leaf/single element is true,
      bAli of a disjoint sum is true if and only if there are exactly two suborders in the disjoint sum,
      one is a chain (bChain = 1) of height 1 or 2 (iHeight1\_2\_More = 1 or 2)
      and the other is a leaf (bLeaf = 1),
      bAli of an order composition is true if and only if the first subtree has bAli = 1
      and all other subtrees have bChain = 1.
\end{itemize}
It may also be nice to compute:
\begin{itemize}
\item an integer value iLongestChain which is the number of elements of a longest chain in the order defined by the modular decomposition up to this node,
      instead of computing iHeight1\_2\_More.
      iLongestChain of a leaf/single element is 1,
      iLongestChain of a disjoint sum is the maximum of the values of iLongestChain for the subtrees of the disjoint sum,
      iLongestChain of an order composition is the sum of the values of iLongestChain for the subtrees of the order composition.
\item a boolean value bAliInverse which is true if and only if the order defined by the modular decomposition up to this node is the inverse/reverse of an ali order.
      bAliInverse of a leaf/single element is true,
      bAliInverse of a disjoint sum is true if and only if there are exactly two suborders in the disjoint sum,
      one is a chain (bChain = 1) of height 1 or 2 (iHeight1\_2\_More = 1 or 2)
      and the other is a leaf (bLeaf = 1),
      bAliInverse of an order composition is true if and only if the last subtree has bAliInverse = 1
      and all other subtrees have bChain = 1.
\item a boolean value bDisjointAli
      which is true if and only if the order defined by the modular decomposition up to this node is a chain of height 1 or 2,
      together with an incomparable element.
      bDisjointAli of a leaf/single element is false, bDisjointAli of a disjoint sum is true if and only if bAli = 1,
      bDisjointAli of an order composition is false.
\item a boolean value bAli11BasedChain
      which is true if and only if the order defined by the modular decomposition up to this node is a (1,1)-based chain.
      bAli11BasedChain of a leaf/single element is false, bAli11BasedChain of a disjoint sum is false,
      bAli11BasedChain of a grouped order composition is true if and only if
        the first subtree has (bDisjointAli = 1 and iHeight1\_2\_More = 1)
      and all other subtrees have bChain = 1,
      bAli11BasedChain of a binary order composition is true if and only if
        the first subtree has (bAli11BasedChain = 1 or (bDisjointAli = 1 and iHeight1\_2\_More = 1))
      and the second subtree has bChain = 1.
\item a boolean value bAli21BasedChain
      which is true if and only if the order defined by the modular decomposition up to this node is a (2,1)-based chain.
      bAli21BasedChain of a leaf/single element is false, bAli21BasedChain of a disjoint sum is false,
      bAli21BasedChain of a grouped order composition is true if and only if
        the first subtree has (bDisjointAli = 1 and iHeight1\_2\_More = 2)
      and all other subtrees have bChain = 1,
      bAli21BasedChain of a binary order composition is true if and only if
        the first subtree has (bAli21BasedChain = 1 or (bDisjointAli = 1 and iHeight1\_2\_More = 2))
      and the second subtree has bChain = 1.
\item a boolean value bAli11EndedChain
      which is true if and only if the order defined by the modular decomposition up to this node is a (1,1)-ended chain
      (the inverse/reverse order of a (1,1)-based chain).
      bAli11EndedChain of a leaf/single element is false, bAli11EndedChain of a disjoint sum is false,
      bAli11EndedChain of a grouped order composition is true if and only if
        the last subtree has (bDisjointAli = 1 and iHeight1\_2\_More = 1)
      and all other subtrees have bChain = 1,
      bAli11EndedChain of a binary order composition is true if and only if
        the second subtree has (bAli11EndedChain = 1 or (bDisjointAli = 1 and iHeight1\_2\_More = 1))
      and the first subtree has bChain = 1.
\item a boolean value bAli21EndedChain
      which is true if and only if the order defined by the modular decomposition up to this node is a (2,1)-ended chain
      (the inverse/reverse order of a (2,1)-based chain).
      bAli21EndedChain of a leaf/single element is false, bAli21EndedChain of a disjoint sum is false,
      bAli21EndedChain of a grouped order composition is true if and only if
        the last subtree has (bDisjointAli = 1 and iHeight1\_2\_More = 2)
      and all other subtrees have bChain = 1,
      bAli21EndedChain of a binary order composition is true if and only if
        the second subtree has (bAli21EndedChain = 1 or (bDisjointAli = 1 and iHeight1\_2\_More = 2))
      and the first subtree has bChain = 1.
\end{itemize}
Computing bChain, iHeight1\_2\_More, bAli/bAliInverse, bDisjointAli, bAli11BasedChain, bAli21BasedChain, bAli11EndedChain,
and bAli21EndedChain on all nodes takes \begin{math}O(n)\end{math} time,
computing iLongestChain on all nodes takes \begin{math}O(n \times \log(n))\end{math} time
(or \begin{math}O?(n)\end{math} time on unit cost RAM-model,
if there is less than \begin{math}2^{64}\end{math} elements which should be the case for efficient computations,
the cost of maximum and sum computation is done by a constant number of hardware instructions on current hardware architectures,
and there will be no empirical asymptotic difference, up to a constant factor, between a \begin{math}O(n)\end{math} and \begin{math}O?(n)\end{math} algorithms
with similar input/output profile (memory access matters a lot)).
)

\section{Orders that are naturally consecutive level-induced suborders}
\label{section:Orders_that_are_naturally_consecutive_level-induced_suborders}

In this section, we assume that a given well-founded order \begin{math}P'\end{math}
is a level-induced suborder of a well-founded order \begin{math}P\end{math}.
We study necessary and sufficient conditions on \begin{math}P'\end{math}
to have that \begin{math}P'\end{math}
is a consecutive level-induced suborder of \begin{math}P\end{math}.
\begin{definition}[nacli orders]
A \emph{nacli order} \begin{math}P'\end{math} is a well-founded order
such that whenever \begin{math}P'\end{math} is isomorphic to a level-induced suborder \begin{math}P''\end{math} of a well-founded order \begin{math}P\end{math},
then \begin{math}P'\end{math} is also isomorphic to a consecutive level-induced suborder of \begin{math}P\end{math}
(or equivalently of \begin{math}P[\InitialSection(P,\Domain(P''))]\end{math}, see the reason for ali orders).
\end{definition}

We first observe that :
\begin{lemma}
For any well-founded order \begin{math}P'\end{math}
containing an induced suborder isomorphic
to \begin{math}O_{obs2} = (\{a,b,c\}, \{a < b, a < c\})\end{math}
or \begin{math}O_{obs3} = (\{a,b,c\}, \{a > b, a > c\})\end{math},
there is a well-founded order \begin{math}P\end{math}
such that \begin{math}P'\end{math} is a level-induced suborder of \begin{math}P\end{math},
but \begin{math}P'\end{math} is not isomorphic to any consecutive level-induced suborder of \begin{math}P\end{math}.
\end{lemma}
\begin{demo}
We use a disjoint sum of well-founded chains of same height to lift each level of \begin{math}P'\end{math}
so that any two levels of \begin{math}P'\end{math} are now \begin{math}\gamma\end{math} levels apart,
where \begin{math}\gamma \geq \omega_{\beta + 1}\end{math},
and the cardinal of \begin{math}\Domain(P')\end{math} is at most \begin{math}\aleph_{\beta}\end{math}.
Again by Zermelo's axiom, there is a bijection \begin{math}f\end{math}
between some ordinal \begin{math}\alpha\end{math} and \begin{math}\Domain(P')\end{math}.
We add a distinct well-founded chain for each element of \begin{math}\Domain(P')\end{math}.
Let \begin{math}\DisjointCopy(i)\end{math} be a chain isomorphic to the ordinal \begin{math}i\end{math}
such that its elements are assumed to be distinct from all other elements considered in the following formula:
\begin{math}\Domain(P) = \Domain(P') \sqcup (\bigsqcup_{i \in \alpha} \DisjointCopy(\gamma \times \Level_{P'}(f(i))))\end{math}.
\begin{itemize}
\item \begin{math}\forall x,y \in \Domain(P'), \OrderFunction(P)(x,y) = \OrderFunction(P')(x,y)\end{math},
\item \begin{math}\forall x,y \in \DisjointCopy(\gamma \times \Level_{P'}(f(i))),
                        \OrderFunction(P)(x,y) = \OrderFunction(\DisjointCopy(\gamma \times \Level_{P'}(f(i))))(x,y)\end{math},
\item \begin{math}\forall x \in \DisjointCopy(\gamma \times \\Level_{P'}(f(i))), \forall y \in \DisjointCopy(\gamma \times \Level_{P'}(f(j))),
                        \OrderFunction(P)(x,y) = \text{ `}\sim\text{'}\end{math},
\item \begin{math}\forall x \in \DisjointCopy(\gamma \times \Level_{P'}(f(i))), \forall y \in \Domain(P'),
        \OrderFunction(P)(x,y) = \text{ `}<\text{' if } i = f^{-1}(y) \text{ or } f(i) < y (\OrderFunction(P')(f(i),y) \in \{=, <\}), \text{ `}\sim\text{' otherwise }\end{math}.
\end{itemize}
Clearly, \begin{math}P'\end{math} is a level-induced suborder of \begin{math}P\end{math},
and any two levels of \begin{math}P'\end{math} are now \begin{math}\gamma\end{math} levels apart,
since \begin{math}\DisjointCopy(\gamma \times \Level_{P'}(f(i)))\end{math} is a longest chain below element \begin{math}f(i)\end{math}.

Let \begin{math}(x,y,z)\end{math} be a triple of elements of \begin{math}\Domain(P)\end{math},
such that \begin{math}x < y, x < z, y \sim z\end{math} (\begin{math}O_{obs2}\end{math}),
or \begin{math}x > y, x > z, y \sim z\end{math} (\begin{math}O_{obs3}\end{math}).
It naturally defines one ordinal \begin{math}\Gap_{P}(x,y,z) = \Sup(\Gap_{P}(x,y), \Gap_{P}(x,z))\end{math}.

Observe that no element of \begin{math}\DisjointCopy(\gamma \times \Level_{P'}(f(i)))\end{math} is more than an element,
unless that element is also in \begin{math}\DisjointCopy(\gamma \times \Level_{P'}(f(i)))\end{math}.
Hence, it cannot be more than two incomparable elements.

Clearly, if it is less than two incomparable elements like \begin{math}x\end{math},
then these two elements are in \begin{math}\Domain(P')\end{math},
and \begin{math}x \in \{f(i)\} \sqcup \DisjointCopy(\gamma \times \Level_{P'}(f(i)))\end{math}
implies that \begin{math}f(i)\end{math} is also less than these two elements.
Moreover, \begin{math}\Gap_{P}(x,y,z) \geq \Gap_{P}(f(i),y,z) \geq \gamma\end{math}, in that case.

Thus, if there is an induced suborder isomorphic
to \begin{math}O_{obs2} = (\{a,b,c\}, \{a < b, a < c\})\end{math} in \begin{math}P'\end{math},
then no consecutive level-induced suborder isomorphic to \begin{math}P'\end{math}
exists in \begin{math}P\end{math},
because of the ordinal gap in \begin{math}P\end{math}
between original levels of \begin{math}P'\end{math}
that is superior to the ordinal corresponding to the cardinal of \begin{math}P'\end{math}.

Otherwise, there is an induced suborder isomorphic to
\begin{math}O_{obs3} = (\{a,b,c\}, \{a > b, a > c\})\end{math}
in \begin{math}P'\end{math}.
Consider such an induced suborder \begin{math}(x,y,z)\end{math} in \begin{math}P\end{math}.
We already noted that \begin{math}x\end{math} must be in \begin{math}P'\end{math};
if both \begin{math}y,z\end{math} are in \begin{math}\DisjointCopy(\gamma \times \Level_{P'}(f(k)))\end{math},
then they are comparable, a contradiction.
Hence, without loss of generality, \begin{math}y \in \{f(j)\} \sqcup \DisjointCopy(\gamma \times \Level_{P'}(f(j)))\end{math},
for some \begin{math}f(j) \neq x, f(j) \in \Domain(P')\end{math}.
It is now trivial to see that \begin{math}(\{x,f(j),z\}, \{x > f(j), x > z\})\end{math}
is also an induced suborder isomorphic to
\begin{math}O_{obs3} = (\{a,b,c\}, \{a > b, a > c\})\end{math}
with \begin{math}\Gap_{P}(x,y,z) \geq \Gap_{P}(x,f(j),z)\end{math}.
But since \begin{math}\Gap_{P}(x,f(j),z) \geq \gamma\end{math},
again we have that no consecutive level-induced suborder isomorphic to \begin{math}P'\end{math}
exists in \begin{math}P\end{math}.
\end{demo}

\begin{corollary}
A nacli order is the disjoint union of well-founded chains.
\end{corollary}

\begin{lemma}
No nacli order has more than one level of size at least 2.
\end{lemma}
\begin{demo}
Again, we create a gap between consecutive levels of \begin{math}P'\end{math}.
We use a unique well-founded chain of height \begin{math}\gamma \times \Height(P')\end{math} to lift all levels of \begin{math}P'\end{math}
so that any two levels of \begin{math}P'\end{math} are now \begin{math}\gamma\end{math} levels apart,
where \begin{math}\gamma \geq \omega_{\beta + 1}\end{math},
and the cardinal of \begin{math}\Domain(P')\end{math} is at most \begin{math}\aleph_{\beta}\end{math}.
\begin{math}\Domain(P) = \Domain(P') \sqcup \DisjointCopy(\gamma \times \Height(P'))\end{math}.
Since added levels have size 1 and original levels are too far apart,
at most one level can have size more than one in a consecutive level-induced suborder.
\end{demo}

\begin{theorem}
A nacli order is a well-founded chain, an antichain, or the disjoint union of a well-founded chain and an antichain.
Equivalently, a nacli order is a \emph{well-founded} order without induced suborder isomorphic
to \begin{math}O_{obs1} = (\{a,b,c,d\}, \{a < b, c < d\}) \equiv \Inv(O_{obs1})\end{math},
\begin{math}O_{obs2} = (\{a,b,c\}, \{a < b, a < c\})\end{math},
or \begin{math}O_{obs3} = (\{a,b,c\}, \{a > b, a > c\})\end{math}.
In particular, nacli orders are a subclass of series parallel interval orders,
and all ali orders except order compositions of a disjoint ali order and a well-founded chain are also nacli orders.
\end{theorem}
\begin{demo}
By previous lemmas, only well-founded chains, antichains, or the disjoint unions of a well-founded chain and an antichain may be nacli orders.
The proof by transfinite induction that such orders are indeed nacli orders is trivial.
In any superorder, fix the first level of the disjoint union of a well-founded chain and an antichain and close the gap with the second level,
then close the gap between the second and third level, etc.
Everything follows from transitivity and the fact that a single well-founded chain can not be lifted by another suborder
that does not contain an isomorphic chain.
\end{demo}
Thus finite nacli orders can be recognized in time \begin{math}O(n+m)\end{math}
with techniques similar to the end of previous section.
(bNacliOfHeight1 = bAntichain is the logical conjunction of bNacliOfHeight1 of subtrees on disjoint sums nodes,
 and false on order compositions nodes.
bNacli is the logical conjunction of bChain on order compositions nodes,
and it is true on disjoint sums nodes if and only if bNacliOfHeight1 is true on all subtrees,
except maybe at most one where instead bChain is true (grouped case),
 or (bChain or bNacli) is true (binary case, can you simplify ``bChain or bNacli''? ;P).)

\section{Algorithms to find ali induced suborders and nacli level-induced suborders}
\label{section:Algorithms_to_find_ali_induced_suborders_and_nacli_level-induced_suborders}

All orders in this section are finite, hence well-founded.
We first start with the simple case of chains and orders made of a chain of height at most 2 and an incomparable element,
i.e. orders that are ali orders and nacli orders at the same time.
Assume we want to find such an order in a superorder \begin{math}P\end{math},
where \begin{math}n\end{math} is the cardinal of \begin{math}\Domain(P)\end{math},
and \begin{math}m\end{math} is the number of comparability relationships in \begin{math}P\end{math}.
We first compute a level decomposition of \begin{math}P\end{math} in time \begin{math}O(n^3)\end{math},
with the additional constraint that we store in each element a reference to another element
that is less than it in the \emph{previous} level.
We do it as follow, once an element has been selected to be added in the current level,
for all elements that are greater than it overwrite their reference with the selected element.
Clearly, the last overwrite will be in the previous level.
It is easy to see that this can be done in time \begin{math}O(m) \leq O(n^2)\end{math}.
Let \begin{math}h\end{math} be the height of \begin{math}P\end{math},
and \begin{math}s\end{math} be the size of the longest chain in the ali and nacli suborder.
\begin{itemize}
\item If the suborder is a chain,
  \begin{itemize}
    \item if \begin{math}h \geq s\end{math}, take any element \begin{math}x\end{math} in the level \begin{math}s - 1\end{math},
    \item otherwise there is no such ((consecutive) level-)induced suborder.
  \end{itemize}
\item Otherwise, for any element \begin{math}x\end{math} in the level \begin{math}l\end{math}
      ranging from \begin{math}s - 1\end{math} to \begin{math}h - 1\end{math} (first loop),
      check if there is an element \begin{math}y\end{math}
      in level \begin{math}l - s + 1\end{math} (second loop) that is incomparable with \begin{math}x\end{math}.
     If no check succeeds, there is no such ((consecutive) level-)induced suborder (this check is sufficient by transitivity).
\end{itemize}
If you got an \begin{math}x\end{math} and optionnaly a corresponding \begin{math}y\end{math},
then you can output the consecutive level-induced suborder made of element \begin{math}y\end{math}
and the chain obtained by following the references set during the level decomposition,
starting from element \begin{math}x\end{math} and iterating \begin{math}s - 1\end{math} times.
Clearly, these two loops take time \begin{math}O(n^2)\end{math}.
Thus, whatever the size of an ali and nacli order,
finding such a ((consecutive) level-)induced suborder has time complexity in \begin{math}O(n^3)\end{math}.

We now look at the odd case of the (1,1)-based chains that are ali orders but not nacli orders.
Assume we want to find such a (1,1)-based chain in a superorder \begin{math}P\end{math}.
We first compute a level decomposition of \begin{math}P\end{math} in time \begin{math}O(n^3)\end{math}.
Then, in time \begin{math}O(m \times \log(n))\end{math}, proceeding from the last level to the first level,
we can compute on each element \begin{math}x\end{math}
the size \begin{math}slc(x)\end{math} of the longest chain starting with \begin{math}x\end{math}:
this size is 1 if no element is greater
and the maximum plus one of \begin{math}slc(y)\end{math} for \begin{math}y\end{math} greater than \begin{math}x\end{math} otherwise,
for backtracking purpose, we keep a reference to such an \begin{math}y\end{math} that gave the maximum for each \begin{math}x\end{math}.
Then, for each element \begin{math}x\end{math} (first loop) such that \begin{math}slc(x)\end{math}
is at least the size of the (1,1)-based chain minus two,
we can enumerate in linear time (second loop) all elements in a level below that are less than \begin{math}x\end{math};
doing this loop on elements in levels below, one level after the another,
we can stop as soon as we find two such elements in the same level,
it takes time \begin{math}O(n^2)\end{math}.
Thus, whatever the size of a (1,1)-based chain,
finding such an induced suborder has time complexity in \begin{math}O(n^3)\end{math}.
The same result applies to (1,1)-ended chains by considering the inverse order.

For (2,1)-based chains that are also ali orders but not nacli orders,
we still use \begin{math}slc(x)\end{math} and loop on elements where it is at least the size of the (2,1)-based chain minus three,
but this time, for the second loop, we search elements \begin{math}y\end{math} less than \begin{math}x\end{math}
and for the third loop, we classify in linear time all elements in a level below the level of \begin{math}y\end{math} into three classes:
\begin{itemize}
\item the elements that are less than \begin{math}y\end{math},
\item the elements that are less than \begin{math}x\end{math} and incomparable with \begin{math}y\end{math},
\item the elements that are incomparable with \begin{math}x\end{math} and \begin{math}y\end{math}.
\end{itemize}
As soon as we have found an element in the first class and an element in the second class in the same level,
we have found a (2,1)-based chain,
and no (2,1)-based chain may exist if we do not find such elements in the same level.
(For the implementation, we just need two element variables elementInClass1, elementInClass2,
initialized with null value,
we loop through the elements one level after the another,
and everytime we change of level, we reinitialize elementInClass1, elementInClass2 with null value.
Whenever the current element is less than \begin{math}y\end{math}, we set elementInClass1 variable;
whenever the current element is less than \begin{math}x\end{math} and incomparable with \begin{math}y\end{math}, we set elementInClass2 variable.
As soon as both variables are distinct of null value, we have a (2,1)-based chain.)
Thus, whatever the size of a (2,1)-based chain,
finding such an induced suborder has time complexity in \begin{math}O(n^3)\end{math}.
The same result applies to (2,1)-ended chains by considering the inverse order.

In order to find (consecutive) level-induced suborders that are nacli orders,
we just modify the algorithm for ali and nacli suborders as follow:
Let \begin{math}r\end{math} be the number of elements of the first level of the suborder minus one.
Replace
\begin{itemize}
\item ``Otherwise, for any element \begin{math}x\end{math} in the level \begin{math}l\end{math}
      ranging from \begin{math}s - 1\end{math} to \begin{math}h - 1\end{math} (first loop),
      check if there is an element \begin{math}y\end{math}
      in level \begin{math}l - s + 1\end{math} (second loop) that is incomparable with \begin{math}x\end{math}.
     If no check succeeds, there is no such ((consecutive) level-)induced suborder (this check is sufficient by transitivity).''
\end{itemize}
by
\begin{itemize}
\item ``Otherwise, for any element \begin{math}x\end{math} in the level \begin{math}l\end{math}
      ranging from \begin{math}s - 1\end{math} to \begin{math}h - 1\end{math} (first loop),
      check if there are \begin{math}r\end{math} elements \begin{math}y_1, \dots, y_r\end{math}
      in level \begin{math}l - s + 1\end{math} (second loop) that are incomparable with \begin{math}x\end{math}.
     If no check succeeds, there is no such (consecutive) level-induced suborder (this check is sufficient by transitivity).''
\end{itemize}

\section{Conclusion}

Maybe we should talk about partial-level-induced suborders,
since we do not impose to keep all elements of a level of a superorder.
However, in that case total/global-level-induced suborders would be rather restricted.
And no similar results could be obtained unless considering superorders of bounded level-width.
\medskip

\noindent H\^atez-vous lentement, et sans perdre courage,\newline
Vingt fois sur le m\'etier remettez votre ouvrage,\newline
Polissez-le sans cesse et le repolissez,\newline
Ajoutez quelquefois, et souvent effacez.

Nicolas Boileau, L'Art po\'etique 1674

\section*{Acknowledgements}
\label{section:acknowledgements}

We thank God: Father, Son, and Holy Spirit. We thank Maria.
They help us through our difficulties in life.

\nocite{*}
\bibliographystyle{abbrvnat}
\bibliography{LL2020OrdreNiveauInduit}
\label{section:bibliography}

\end{document}